\documentclass{amsart}
\usepackage{amsthm, amssymb, amsmath, amscd, amsfonts}
\usepackage{xypic, graphicx, epsf, mathptmx, verbatim}
\title[String topology prospectra and Hochschild cohomology]{String topology prospectra and Hochschild cohomology}
\author[Kate Gruher and Craig Westerland]{Kate Gruher and Craig Westerland}
\input xy

\newtheorem{theorem}{Theorem}[section]
\newtheorem{proposition}[theorem]{Proposition}

\newtheorem{lemma}[theorem]{Lemma}

\newtheorem{corollary}[theorem]{Corollary}

\theoremstyle{definition}
\newtheorem{definition}[theorem]{Definition}

\theoremstyle{remark}
\newtheorem{remark}[theorem]{Remark}
\newtheorem{example}[theorem]{Example}

\def\R{\ensuremath{\mathbb{R}}}

\def\N{\ensuremath{\mathbb{N}}}

\def\L{\ensuremath{\mathbb{L}}}
\def\g{\ensuremath{\mathfrak{g}}}

\def\id{{\rm id}}

\def\Hom{{\rm Hom}}
\def\RHom{{\rm RHom}}

\def\Tot{{\rm Tot}}

\def\smash{\land}
\def\scs{\scriptstyle}
\def\into{\hookrightarrow}

\def\colim{\ensuremath{{\rm colim} \,}}

\begin{document}
\bibliographystyle{amsalpha}
\maketitle
\begin{abstract}
We study string topology for classifying spaces of connected compact Lie groups, drawing connections with Hochschild cohomology and equivariant homotopy theory.  First, for a compact Lie group $G$, we show that the string topology prospectrum $LBG^{-TBG}$ is equivalent to the homotopy fixed-point prospectrum for the conjugation action of $G$ on itself, $G^{hG}$.  Dually, we identify $LBG^{-ad}$ with the homotopy orbit spectrum $(DG)_{hG}$, and study ring and co-ring structures on these spectra.  Finally, we show that in homology, these products may be identified with the Gerstenhaber cup product in the Hochschild cohomology of $C^*(BG)$ and $C_*(G)$, respectively.  These, in turn, are isomorphic via Koszul duality.
\end{abstract}
%
%
%
%
%
\section{Introduction} \label{intro}

Let $G$ be a connected compact Lie group.  The free loop space
$$LBG := Map(S^1, BG)$$
of the classifying space of $G$ is a natural object of study for topologists, representation theorists, and mathematical physicists.  Its $K$-theory is related to an important example of a topological field theory, the Verlinde algebra of positive energy representations of the loop group $LG$ \cite{fht3}.  In this article we study $LBG$ and natural field-theoretic algebraic structures which it supports from several points of view -- string topology, Hochschild cohomology, and equivariant stable homotopy theory.


\subsection{Equivalences of (pro-)spectra}

In string topology, one studies the free loop space $LM$ of a closed, oriented, finite dimensional manifold $M$.  Using a combination of intersection theory on $M$ and concatenation of loops with common basepoints, Chas and Sullivan \cite{cs} gave the shifted homology of $LM$ the structure of a Gerstenhaber algebra.  The ring structure was reinterpreted in the language of stable homotopy theory by Cohen and Jones in \cite{cj} in the form of a (Thom) ring spectrum $LM^{-TM}$. 

Although $BG$ is not a finite dimensional manifold, it does admit a filtration by finite dimensional manifolds.  In \cite{gs}, Salvatore and the first author defined an inverse system of ring spectra (or pro-ring spectrum) $LBG^{-TBG}$ using this filtration and analogues of the string topology techniques of \cite{cs, cj}.  In \cite{goperads} the second author studied a ring spectrum $G^{hG}$, the homotopy fixed point spectrum for the action of $G$ on itself by conjugation.  This spectrum is best understood as a pro-ring spectrum.  One purpose of this paper is to show that there is an equivalence between the geometrically constructed $LBG^{-TBG}$ and $G^{hG}$, whose description is equivariant stable homotopy-theoretic.

\begin{theorem} \label{equiv_thm}

The transfer map $\tau^G$ defines an equivalence of pro-ring spectra
$$LBG^{-TBG} \simeq G^{hG}.$$

\end{theorem}

One should compare this result to \cite{klein_fiber}, where Klein shows that for a Poincar\'e duality group $G$ with classifying space $M = BG$ a Poincar\'e duality space of formal dimension $d$, there is an equivalence of the \emph{spectrum} $G^{hG}$ (not a pro-spectrum) with the string topology spectrum $LM^{-TM}$.

It is worth pointing out that the spectrum $G^{hG}$ is equivalent to $THH^\bullet (G, G)$, the topological Hochschild cohomology of (the suspension spectrum of) $G$.  This foreshadows Theorem \ref{hoch_thm} below.

In \cite{kate}, the first author showed that the prospectrum
$LBG^{-TBG}$ is Spanier-Whitehead dual (in the sense of Christensen
and Isaksen \cite{ci}) to a spectrum $LBG^{-ad}$.  There is a
coproduct on that spectrum which, upon application of a cohomology
theory, gives an (untwisted) analogue of the Freed-Hopkins-Teleman product in twisted equivariant $K$-theory (or fusion product in the Verlinde algebra).

In light of this duality and Theorem \ref{equiv_thm}, the following should be unsurprising.

\begin{theorem} \label{dual_orbit_thm}

There is an equivalence of co-ring spectra $LBG^{-ad} \simeq (DG)_{hG}$.

\end{theorem}

Here $DG = F(\Sigma^{\infty} G_+, S^0)$ is the Spanier-Whitehead dual of $G$, equipped with a naive $G$-action dual to the conjugation action on $G$.  We describe the coproduct on the Borel construction $DG_{hG} = EG_+ \smash_G DG$ in section \ref{coproduct_section} below.

A remark on terminology is in order.  Throughout this paper, the terms ``ring spectrum" and ``pro-ring spectrum" will be used to describe objects whose multiplication is associative up to homotopy.  For more highly structured ring spectra, we will employ the $S$-algebras of \cite{EKMM}.  Additionally, the term ``pro-ring spectrum" (resp. ``pro-$S$-algebra") denotes an inverse system of ring spectra (resp. $S$-algebras), rather than a monoid in the category of prospectra.

Further, we will not consider strict co-ring spectra, and only require them to be co-associative up to homotopy.  Indeed, for most of this paper, we work in the homotopy category.  However, the prospectrum $G^{hG}$ is a (strict) pro-$S$-algebra, so Theorem \ref{equiv_thm} can be thought of as a rectification result for $LBG^{-TBG}$.  This answers in the affirmative Conjecture 10 of \cite{kate}.

\subsection{Homological computations}

A natural question is how to compute the (co)homology of these
(pro)spectra..  Let $k$ be a field; all of our (co)chain and (co)homology groups will have coefficients in $k$. 

Our approach is through Hochschild cohomology.  For a
differential graded algebra $A$ and a dg $A$-module $M$, $HH_*(A, M)$
and $HH^*(A, M)$ are the Hochschild homology and cohomology of $A$
with coefficients in $M$.  Recall that for any topological group $K$
and topological space $X$, there are isomorphisms
$$\begin{array}{ccc}
H_*(LBK) \cong HH_*(C_*(K), C_*(K))& {\rm and} & H^*(LX) \cong
HH_*(C^*(X), C^*(X))
\end{array}$$
where $C_*(K)$ is given the structure of a dga via the Pontrjagin
product, and $C^*(X)$ via the cup product of cochains.

In \cite{cj}, Cohen-Jones modified the latter isomorphism to give an
isomorphism of rings
$$H_*(LM^{-TM}) \cong HH^*(C^*(M), C^*(M))$$
for finite dimensional manifolds $M$.  We adapt both of these
computations to the context of string topology on $BG$.

\begin{theorem} \label{hoch_thm}

If $G$ is a connected compact Lie group, the following rings are
mutually isomorphic:

\begin{enumerate}

\item $H_*^{pro}(LBG^{-TBG})$, with the string topology product of
  \cite{gs}. \label{gs_item}

\item $H^{-*}(LBG^{-ad})$, with the ring structure induced by the ``fusion''
  coproduct on $LBG^{-ad}$, defined in \cite{kate}. \label{fusion_item}

\item $HH^*(C^*(BG), C^*(BG))$, with the Gerstenhaber cup product. \label{hoch_BG_item}

\item $HH^*(C_*(G), C_*(G))$, with the Gerstenhaber cup product. \label{hoch_G_item}

\end{enumerate}

\end{theorem}

In (\ref{gs_item}), $H_*^{pro}(LBG^{-TBG})$ denote the inverse limit of the homologies of the terms in the prospectrum $LBG^{-TBG}$.

Here is a summary of the proof.  To show the equivalence of (\ref{gs_item}) and (\ref{fusion_item}), one uses the Spanier-Whitehead duality result of \cite{kate}.  The isomorphism of
the rings in (\ref{gs_item}) and (\ref{hoch_BG_item}) uses, as in
\cite{cj}, a cosimplicial model for $LBG^{-TBG}$.  Finally, the differential graded algebras $C^*(BG)$ and $C_*(G)$ are
Koszul (or cobar) dual: $C^*(BG)$ is equivalent to the cobar complex
for the differential graded algebra $C_*(G)$ (and vice versa).  As
Hochschild cohomology is insensitive to Koszul duality \cite{fmt, po}, one
obtains an isomorphism of the rings in (\ref{hoch_BG_item}) and (\ref{hoch_G_item}).

Write this collection of isomorphisms in the following form:
$$\xymatrix{
H_*^{pro}(LBG^{-TBG}) \ar@{<->}[r] \ar@{<->}[d]_D & HH^*(C^*(BG),
C^*(BG)) \ar@{<->}[d]^-K \\
H^{-*}(LBG^{-ad}) \ar@{<->}[r] & HH^*(C_*(G), C_*(G))
}$$
In this diagram, the horizontal isomorphisms are ``geometric'' in the
sense that they come from explicit models for the spectra involved.  The vertical isomorphism $D$ is induced by Spanier-Whitehead duality,
and $K$ is induced by Koszul duality.  Consequently one may interpret this theorem as saying that the Spanier-Whitehead duality (of \cite{kate}) between the Chas-Sullivan and Freed-Hopkins-Teleman products is manifested in Hochschild cohomology as an aspect of Koszul duality.

Recent work of Vaintrob \cite{vaintrob} gives an analogue of the
isomorphisms between (\ref{gs_item}) and (\ref{hoch_G_item}) in the related case that $M^n=BG$ is a closed, oriented, aspherical manifold, and $G = \pi_1(M)$ is a discrete group.  Namely, Vaintrob gives an isomorphism of BV algebras
$$HH^{n-*} (k[\pi_1(M)], k[\pi_1(M)]) \cong \Sigma^{-n} H_*(LM)$$
where $k[\pi_1(M)] \cong C_*(G; k)$ is the group algebra on the fundamental group of $M$.

Similar multiplicative structures coming from Chen-Ruan cohomology and string topology of orbifolds and stacks have been studied recently (see, e.g., \cite{glsux, bgnxinertia}).  In the final part of this paper, we relate these constructions to the algebras described above.

We would like to thank Paul Bressler, Ralph Cohen, and Jesper Grodal for stimulating conversations on this material.

\section{The pro-ring spectra}

Let us review the construction of these pro-ring spectra.  Both will be defined 
using a filtration of $EG$ -- a contractible space upon which $G$ acts freely -- 
by finite dimensional free $G$-manifolds.  To do this, we proceed as follows.  
Because $G$ is compact Lie, there exists a finite-dimensional, faithful 
representation $V$ of $G$.

\begin{definition}

Define $EG_n$ to be the space of linear embeddings of $V$ into $\R^n$.

\end{definition}

Since the action of $G$ on $V$ is faithful, when $EG_n$ is nonempty it is a free 
$G$-space, so fits into a principal $G$-fibration
$$G \to EG_n \to BG_n,$$
where we define $BG_n:= EG_n/G$.  Furthermore, by definition, $EG_n$ and $BG_n$ 
are both smooth manifolds.  Finally, the filtered union of the sequence
$$EG_1 \subseteq \dots \subseteq EG_n \subseteq EG_{n+1} \subseteq \dots$$
is the space of linear embeddings of $V$ into $R^{\infty}$ and contractible, so 
is therefore a model for $EG$; i.e.
$$\colim EG_n  = EG.$$
Similarly, $\colim BG_n = BG$.

\begin{example}

For instance, when $G=SO(k)$, $V$ may be taken to be $\R^k$, with the defining 
action of $SO(k)$ on $V$.  Then $BG_n$ is the Grassmannian of $k$-planes in 
$\R^n$, and $EG_n$ is the corresponding Stiefel manifold.

\end{example}

\subsection{The string topology of $BG$}

\begin{definition}

Let $Ad(EG_n)$ denote the total space of the principal $G$-bundle
$$\pi: EG_n \times_G G \to BG_n,$$
where $G$ acts on itself by conjugation.

\end{definition}

Since $BG_n$ is a manifold, it has a tangent bundle, which one can pull back to 
$Ad(EG_n)$ via $\pi$.  In \cite{gs}, it was shown that the Thom spectra
$$Ad(EG_n)^{-TBG_n} := Ad(EG_n)^{-\pi^*(TBG_n)}$$
are ring spectra, using a construction analogous to
Cohen-Jones' construction of string topology operations in \cite{cj}.
Specifically, one has a commutative diagram:
$$\xymatrix{
G \times G \ar[d] & G\times G \ar[l]_-{=} \ar[r]^-{\mu} \ar[d] & G
\ar[d] \\
Ad(EG_n) \times Ad(EG_n) \ar[d] & Ad(EG_n) \times_{BG_n} Ad(EG_n)
\ar[l]_-{\tilde{\Delta}} \ar[r]^-{\tilde{\mu}} \ar[d] & Ad(EG_n) \ar[d]
\\
BG_n \times BG_n & BG_n \ar[l]_-{\Delta} \ar[r]^-{=} & BG_n
}$$
Because $\Delta$ is finite codimension, so too is $\tilde{\Delta}$;
hence both admit \emph{umkehr} Pontrjagin-Thom collapse maps.  Multiplication in the spectrum $Ad(EG_n)^{-TBG_n}$ is given by the
composite of the Pontrjagin-Thom collapse for $\tilde{\Delta}$ with $\tilde{\mu}$.

Furthermore, the natural inclusions $EG_n \subseteq EG_{n+1}$ define (via associated Pontrjagin-Thom maps) a tower of ring spectra
$$Ad(EG_1)^{-TBG_1} \gets \dots \gets Ad(EG_n)^{-TBG_n} \gets 
Ad(EG_{n+1})^{-TBG_{n+1}} \gets \dots .$$
Since there is a homotopy equivalence
$$LBG \simeq Ad(EG),$$
this pro-ring spectrum is denoted $LBG^{-TBG}$.

\subsection{The naive homotopy fixed point prospectrum}\label{naivefps}

Using the manifolds $EG_n$, one can define another pro-ring spectrum.  Consider the 
function spectrum
$$F({EG_n}_+, \Sigma^{\infty} G_+)^G$$
of $G$-equivariant maps from $EG_n$ to the suspension spectrum of $G$.  Here 
$\Sigma^{\infty} G_+$ is regarded as a naive $G$-spectrum, with conjugation 
action.  This may be given a ring product $\mu_n$ using the following diagram:
$$\xymatrix{
F({EG_n}_+, \Sigma^{\infty} G_+)^G \smash F({EG_n}_+, \Sigma^{\infty} G_+)^G 
\ar[dd]_{\mu_n} \ar[r]^-{smash} & F({EG_n \times EG_n}_+, \Sigma^{\infty} 
{G\times G}_+)^{G \times G} \ar[d]^{i} \\
 & F({EG_n \times EG_n}_+, \Sigma^{\infty} {G\times G}_+)^{G} \ar[d]^{\Delta^*} 
\\
F({EG_n}_+, \Sigma^{\infty} G_+)^G & F({EG_n}_+, \Sigma^{\infty} {G\times 
G}_+)^{G} \ar[l]^-{\mu_*} \\
}$$
Here $smash$ smashes two functions together.  The spectrum 
$$F({EG_n \times EG_n}_+, \Sigma^{\infty} {G\times G}_+)^{G}$$
is the space of maps that are equivariant with respect to the diagonal $G$ 
action on each factor, so $i$ is a forgetful map.  The diagonal 
$$\Delta: EG_n \to EG_n \times EG_n$$
is a $G$-equivariant map, so induces $\Delta^*$.  Similarly, $\mu_*$ is induced 
by the multiplication $\mu:G \times G \to G$, which is a $G$-equivariant 
map (with respect to the diagonal action by conjugation).

It was shown in \cite{goperads} that $\mu_n$ makes $F({EG_n}_+, \Sigma^{\infty} 
G_+)^G$ into an associative $S$-algebra (in fact, it is the first term of an 
operad in the stable category).

The natural inclusions $EG_n \subseteq EG_{n+1}$ are $G$-equivariant, so induce 
maps of $S$-algebras
$$F({EG_n}_+, \Sigma^{\infty} G_+)^G \gets F({EG_{n+1}}_+, \Sigma^{\infty} 
G_+)^G$$
which assemble into the pro-$S$-algebra
$$F({EG_1}_+, \Sigma^{\infty} G_+)^G \gets \dots \gets F({EG_n}_+, 
\Sigma^{\infty} G_+)^G \gets F({EG_{n+1}}_+, \Sigma^{\infty} G_+)^G 
\gets \dots .$$
For a naive $G$-spectrum $X$, the function spectrum $F(EG_+, X)^G$ is called the \emph{homotopy fixed point spectrum} $X^{hG}$.  We will therefore denote this pro-$S$-algebra $\Sigma^{\infty} G_+^{hG}$.  For brevity, we will tend to refer to it simply as $G^{hG}$.

\subsection{An alternate homotopy fixed point prospectrum}

In equivariant stable homotopy theory there is another notion of suspension 
spectrum.  For a space $X$, one may define the spectrum $\Sigma^{\infty}_G X$ 
whose $n^{\rm th}$ space is 
$$Q_G \Sigma^n X = \colim_V \Omega^V \Sigma^{V+\R^n} X,$$
and the colimit is taken over a complete $G$-universe of real finite-dimensional representations $V$ of $G$.  Here $S^{V} = V \cup \{ \infty \}$ is the one-point compactification of $V$,
$$\Sigma^{V+\R^n} X = S^V \smash S^n \smash X$$
and $\Omega^V Y = F(S^V, Y)$ is the function space of based continuous maps $S^V \to Y$.  We make $\Sigma^{\infty}_G X$ into a naive $G$-spectrum as follows: for $f \in \Omega^V \Sigma^{V+\R^n} X$, $g \in G$, and $v \in S^{V} = V \cup \{ \infty \}$,
$$(g \cdot f)(v) = g f(v \cdot g^{-1}).$$
This extends over the colimit to give an action on each term of the spectrum.  
This, in turn, assembles into a naive action of $G$ on $\Sigma^{\infty}_G X$.

Replacing $\Sigma^{\infty} G_+$ with $\Sigma^{\infty}_G G_+$ above (and using 
precisely the same arguments), we get a pro-$S$-algebra
$$F({EG_1}_+, \Sigma^{\infty}_G G_+)^G \gets \dots \gets F({EG_n}_+, 
\Sigma^{\infty}_G G_+)^G \gets F({EG_{n+1}}_+, \Sigma^{\infty}_G G_+)^G 
\gets \dots .$$
We will denote this pro-$S$-algebra $\Sigma^{\infty}_G G_+^{hG}$.

There is a natural map
$$e: \Sigma^{\infty} X \to \Sigma^{\infty}_G X,$$
for one can regard the terms of $\Sigma^{\infty} X$ as a similar colimit, only 
taken over the family of a trivial $G$-representations.  This map is an 
equivariant map which is a nonequivariant equivalence \cite{acd, gm}
and thus gives 
an equivalence on homotopy fixed points.  Consequently the induced map of 
prospectra
$$e: \Sigma^{\infty} G_+^{hG} \to \Sigma^{\infty}_G G_+^{hG}$$
is a pro-equivalence.

\section{The co-ring spectra} \label{coproduct_section}

In this section we study the spectra $LBG^{-ad}$ and $(DG)_{hG}$ and
the coproducts defined on each.

\subsection{The spectrum $LBG^{-ad}$}

We recall the definition of $LBG^{-ad}$.  Let $\g$ be the Lie
algebra of $G$, equipped with the adjoint action of $G$.  Then one may
form a flat bundle $ad$ over $Ad(EG) = EG \times_G G$ with total space
$$ad := (EG \times G \times \g)/G$$
The Thom spectrum of the virtual bundle $-ad$ over $Ad(EG) \simeq LBG$
is what we shall call $LBG^{-ad}$.
%

\subsection{Group actions on variants of $DG$}

The group action of $G$ on itself by conjugation induces a naive action of $G$ on $DG = F(\Sigma^{\infty} G_+, S^0)$ by pre-conjugation.  We explore two variants on this action that are more geometrically defined.

The tangent bundle $TG$ of $G$ can be given the structure of a $G$-equivariant vector bundle, lifting the conjugation action on $G$: for $g \in G$ define $c_g : G \to G$ to be conjugation by $g$.  For $h \in G$ and $v \in T_hG$, we define 
$$g \cdot (h, v) := (c_g(h), d_h(c_g)(v))$$
where $d_h(c_g)$ is the derivative of $c_g$ at $h$.  This construction makes the Thom spectrum $G^{-TG}$ into a naive $G$-spectrum.

Alternatively, consider the Lie algebra $\g := T_eG$ alone.  It inherits an action of $G$ as a subspace of $TG$; this is the adjoint action.  This makes $S^{\g} = \g \cup \{ \infty \}$ into a $G$-space, and thus $S^{-\g}$ a naive $G$-spectrum.  Smashing with the conjugation action on $G$ gives a naive action of $G$ on $S^{-\g} \smash G_+$.

\begin{proposition} \label{equiv_prop}

There are equivariant equivalences
$$DG \simeq G^{-TG} \simeq S^{-\g} \smash G_+$$

\end{proposition}

\begin{proof}

The first equivalence is Atiyah duality.  The second follows from the fact that $G$ is parallelizable.

\end{proof}

Notice that this gives an alternate construction of $LBG^{-ad}$; from
the construction of $ad$, it is apparent that
$$(S^{-\g} \smash G_+)_{hG} = LBG^{-ad}$$
Then Proposition \ref{equiv_prop} implies part of Theorem \ref{dual_orbit_thm}: taking homotopy orbits we see that
$$(DG)_{hG} \simeq (S^{-\g} \smash G_+)_{hG} = LBG^{-ad}$$

\subsection{Co-ring spectra}

In \cite{kate}, it was shown that $Ad(E)^{-ad}$ is a co-ring spectrum,
when $p: E \to M$ is a principal $G$-bundle over a finite dimensional
manifold $M$.  It is not hard to extend this to the
infinite-dimensional case $M=BG$:

\begin{proposition} \label{coring_prop}

The spectrum $LBG^{-ad} = (S^{-\g} \smash G_+)_{hG}$ is a homotopy co-associative co-ring spectrum.

\end{proposition}

\begin{proof}

The multiplication map $m: G \times G \to G$ is a principal $G$-bundle; the fibre over the identity is $\{ (g, g^{-1}), g \in G \} \cong G$.  Consequently there is a (stable) transfer map
$$m_!: S^\g \smash G_+ \to G \times G_+$$
which is well-defined up to homotopy.  If we give $G \times G$ an action of $G$ by conjugation in each factor it is easy to see that $m$ is equivariant.  Therefore $m_!$ is also.  Smashing with $S^{-2 \g}$ and taking homotopy orbits gives
$$M: (S^{-\g} \smash G_+)_{hG} \to ((S^{-\g} \smash G_+) \smash (S^{-\g} \smash G_+))_{hG}$$
Here $M = (\id_{S^{-2 \g}} \smash m_!)_{hG}$.

For any two naive $G$-spectra $X$ and $Y$, there is a natural map
$$d:(X \smash Y)_{hG} \to X_{hG} \smash Y_{hG}$$
induced by the diagonal on $EG$.  We may define the coproduct on $LBG^{-ad}$ to be the composite $d \circ M$.

To see that the coproduct is co-associative, first observe that $(m_! \wedge \id) \circ m_! = (\id \wedge m_!) \circ m_!$ as maps
$$S^\g \smash S^\g \smash G_+ \to S^\g \smash G_+ \smash G_+ \to G_+ \smash G_+ \smash G_+ $$
since both are equal to the transfer map for the principal $G \times G$-bundle $G\times G\times G \to G$ given by three-term multiplication.  Smashing with $S^{-3\g}$ and taking homotopy orbits shows that the two compositions in the diagram below are equal.
$$\xymatrix{
& ((S^{-\g}\smash G_+) \smash (S^{-\g} \smash G_+) \smash (S^{-\g}\smash G_+))_{hG} \\
(S^{-\g}\smash G_+)_{hG} \ar[r]^-{M} & ((S^{-\g}\smash G_+) \smash (S^{-\g}\smash G_+))_{hG} \ar[u]_{ m_!\smash \id} \ar[d]^{\id \smash m_!} \\
& ((S^{-\g}\smash G_+) \smash (S^{-\g} \smash G_+) \smash (S^{-\g}\smash G_+))_{hG}
}$$
Co-associativity then follows from the naturality of $d$ and the co-associativity of the diagonal map on $EG.$

\end{proof}

Since $G$ is a finite complex, the Spanier-Whitehead dual $DG$ is also equipped with a natural co-ring spectrum structure, dual to the multiplication $m$ in $G$.  Since $m$ is $G$-equivariant (with respect to the diagonal conjugation action), the coproduct on $DG$ is also equivariant.  This allows us to define a coproduct on the Borel construction $(DG)_{hG}$ by
$$\xymatrix@1{(DG)_{hG} \ar[r]^-{Dm} & (DG \smash DG)_{hG} \ar[r]^-{d} & (DG)_{hG} \smash (DG)_{hG}}$$
It is evident that the Atiyah-duality equivalence $(DG)_{hG} \simeq (S^{-\g} \smash G_+)_{hG} = LBG^{-ad}$ of Proposition \ref{equiv_prop} respects these co-ring structures.  This completes the proof of Theorem \ref{dual_orbit_thm}.

\section{The proof of Theorem \ref{equiv_thm}}

We begin the proof of Theorem \ref{equiv_thm} with the following lemma which 
asserts that the terms in each prospectrum are equivalent.

\begin{lemma}

The transfer map for the principal $G$-fibration
$$p: EG_n \times G \to Ad(EG_n)$$
gives rise to an equivalence
$$\tau_n: Ad(EG_n)^{-TBG_n} \simeq F({EG_n}_+, \Sigma^{\infty}_G G_+)^G.$$

\end{lemma}

\begin{proof}

Write $\g$ for the Lie algebra of $G$ and give it the adjoint $G$-action.  Then 
one may form the vector bundle
$$(EG_n \times G \times \g)/G$$ 
over $Ad(EG_n)$, with fibre $\g$.  We will write the Thom space of this bundle 
as $Ad(EG_n)^{\g}$.

Recall from \cite{madsch} that the transfer map $\tau^G$ is an equivalence of 
spectra 
$$\tau^G: \Sigma^{\infty} Ad(EG_n)^{\g} \to (\Sigma^{\infty}_G EG_n 
\times G_+)^G.$$
Let $T$ denote the tangent bundle of $Ad(EG_n)$, and $p^*(T)$ its pullback to 
$EG_n \times G$ via $p$.  Then $\tau^G$ extends to an equivalence of Thom spectra
$$\tau_n = (\tau^G)^{-T}: Ad(EG_n)^{\g-T} \to (\Sigma^{\infty}_G (EG_n 
\times G)^{-p^*(T)})^G. \leqno{(*)}$$

There is a splitting of the tangent bundle of $Ad(EG_n)$:
$$T = \g \oplus \pi^* (TBG_n).$$
$BG_n$ embeds as the unit section of the projection $\pi$, and the vertical 
tangent bundle to $\pi$ is $\g$.  Therefore the lefthand side of $(*)$ may be 
written as $Ad(EG_n)^{-TBG_n}$.

Examine the righthand side of $(*)$.  Since $p$ is a $G$-principal fibration, we 
know that
$$p^*(T) \oplus \g  = T(EG_n \times G),$$
and here $\g$ is a trivial bundle over $EG_n \times G$.  Note that $\g$ is the 
lift of the tangent bundle of $G$ to $EG_n \times G$.  Therefore $p^*(T)$ is 
stably equivalent to the lift of $TEG_n$ to $EG_n \times G$.  Therefore the 
righthand side of $(*)$ may be written as
$$(EG_n^{-TEG_n} \smash \Sigma^{\infty}_G G_+)^G.$$

Atiyah duality then tells us that, since $EG_n$ is a manifold, $EG_n^{-TEG_n}$ 
is the Spanier-Whitehead dual of ${EG_n}_+$:
$$EG_n^{-TEG_n} \simeq F(\Sigma^{\infty} {EG_n}_+, S^0).$$
Using this along with the fact that for a finite spectrum $X$, there is an 
equivalence
$$F(Y, X) \simeq F(Y, S^0) \smash X,$$
we see that the righthand side of $(*)$ is
$$F({EG_n}_+, \Sigma^{\infty}_G G_+)^G.$$

\end{proof}

\begin{lemma}

The maps $\tau_n$ are maps of ring spectra, up to homotopy.

\end{lemma}

\begin{proof}

We will show that the following diagram homotopy commutes.
%
%
%
%
%
$$\xymatrix{
\scs Ad(EG_n)^{-TBG_n} \land Ad(EG_n)^{-TBG_n}  \ar[r]^-{\scs \tau_n \land \tau_n}  \ar[d]_{\tau_G^{G\times G}} & \scs (EG_n^{-TEG_n} \land \Sigma^{\infty}_G G_+)^G \land (EG_n^{-TEG_n} \land \Sigma^{\infty}_G G_+)^G  \ar[d]^{j \circ i \circ smash} \\
\scs Ad_G(EG_n \times EG_n)^{-(\g \oplus T(BG_n\times BG_n))}  \ar[r]^-{\scs \tau_n '}  \ar[d]_ {\tau_{\Delta'}} & \scs ((EG_n \times EG_n)^{-T(EG_n \times EG_n)} \land \Sigma^{\infty}_G G\times G_+)^G  \ar[d]^{\Delta^*} \\
\scs (Ad(EG_n) \times _{BG_n} Ad(EG_n))^{-TBG_n}  \ar[r]^-{\scs \tau_n ''}  \ar[d]_{\mu_*} & \scs (EG_n^{-TEG_n} \land \Sigma^{\infty}_G G\times G_+)^G \ar[d]^{\mu_*} \\
\scs Ad(EG_n)^{-TBG_n} \ar[r]^-{\tau_n} & \scs (EG_n^{TEG_n} \land \Sigma^{\infty}_G G_+)^G
}$$
Here, 
$$Ad_G(EG_n \times EG_n) = (EG_n \times G \times EG_n \times G)/G$$
where $G$ acts diagonally.  All of the horizontal maps are transfer maps:
$\tau_n'$ is the transfer for the principal $G$-bundle
$$EG_n\times G \times EG_n \times G \to (EG_n \times G \times EG_n \times G)/G,$$
Thomified with respect to the bundle $-(T \times T)$, and $\tau_n''$ is the transfer for 
$$EG_n\times G \times G \to (EG_n\times G \times G)/G = Ad(EG_n)\times _{BG_n} Ad(EG_n),$$
Thomified with respect to $-(\g \oplus \pi^*(TBG_n)).$

First consider the top square. 
The map $\tau_G^{G\times G}$ is a transfer map similar to $\tau^G$, 
arising from a Pontrjagin-Thom collapse map.  
Here is it Thomified with respect to $-(TBG_n \times TBG_n)$.  
The map $j$ is induced by the natural map 
$$\Sigma_G^{\infty} G_+ \land \Sigma_G^{\infty} G_+ \to \Sigma_G^{\infty}G\times G_+.$$
This top square commutes by the subgroup naturality of the transfer construction \cite{madsch}.

Next, consider the middle square. 
The map $\Delta^*$ is the Spanier-Whitehead dual of the diagonal 
$\Delta: EG_n \into EG_n \times EG_n,$ 
hence is the Pontrjagin-Thom collapse map for the embedding $\Delta$.  
Likewise, $\tau_{\Delta'}$ is the Pontrjagin-Thom collapse map for 
$$\Delta': (EG_n \times G \times G)/G \into (EG_n \times G \times EG_n \times G)/G.$$
The transfer maps are also collapse maps, and the two ways around this square are the same collapse map, up to homotopy.

In the third square, both vertical maps are induced by the group multiplication on $G$.  Since this multiplication is equivariant for the diagonal conjugation action, the bottom square commutes \cite{madsch}.

Notice that the composition $\tau_{\Delta'} \circ \tau_G^{G\times G}$ is the Pontrjagin-Thom collapse map for the embedding
$$\tilde{\Delta}: Ad(EG_n) \times_{BG_n} Ad(EG_n) \into Ad(EG_n) \times Ad(EG_n).$$
Thus, the composition 
$\mu_* \circ \tau_{\Delta'} \circ \tau_G^{G\times G}$ 
is the same as the ring spectrum multiplication on $Ad(EG_n)^{-TBG_n}$ given in \cite{gs}.  
Furthermore, after identifying 
$$(EG_n^{-TEG_n} \land \Sigma_G^{\infty}G_+)^G \simeq F({EG_n}_+, \Sigma_G^{\infty}G_+)^G,$$
we see that the product given by 
$\mu_* \circ\Delta^*\circ j\circ i \circ smash$ 
is the same as that defined in section~\ref{naivefps}.  
Thus $\tau_n$ is a map of ring spectra, up to homotopy.

\end{proof}

\begin{lemma}

The maps $\tau_n$ commute with the maps defining the prospectra $LBG^{-TBG}$ and 
$G^{hG}$.  That is, they define a map of prospectra.

\end{lemma}

\begin{proof}

First observe that the structure maps 
$$F({EG_n}_+, \Sigma_G^{\infty} G _+)^G \gets F({EG_{n+1}}_+, \Sigma_G^{\infty} G _+)^G $$
define maps
$$(EG_n^{-TEG_n} \land \Sigma_G^{\infty}G_+)^G \gets (EG_{n+1}^{-TEG_{n+1}} \land \Sigma_G^{\infty}G_+)^G$$
which are induced by the Spanier-Whitehead dual of the inclusions $EG_n \subseteq EG_{n+1}$ and hence are the corresponsing Pontrjagin-Thom collapse maps.
We need to check that the diagram
\begin{equation}\label{collapses}
\xymatrix{
Ad(EG_n)^{-TBG_n}  \ar[r]^-{\tau_n} & (EG_n^{-TEG_n}\land \Sigma_G^{\infty}G_+)^G 
\\
Ad(EG_{n+1})^{-TBG_{n+1}}  \ar[r]^-{\tau_{n+1}} \ar[u] & (EG_{n+1}^{-TEG_{n+1}}\land \Sigma_G^{\infty}G_+)^G \ar[u]
}
 \end{equation}
commutes.  From the construction of the transfer map we have a commutative diagram
$$\xymatrix{
Ad(EG_n)^{\g \oplus \nu} \ar[r]^-{\tau_n} & (EG_n^{\nu}\land \Sigma_G^{\infty}G+1)^G
\\
Ad(EG_{n+1})^{\g} \ar[r]^-{\tau_{n+1}} \ar[u] & ({EG_{n+1}}_+ \land \Sigma_G^{\infty}G_+)^G \ar[u]
}$$
where the vertical maps are the collapse maps and $\nu$ is the pullback of the normal bundle of $BG_n$ in $BG_{n+1}.$  Thomifying the diagram above with respect to $-T(Ad(EG_{n+1}))$ yields the diagram (\ref{collapses}).

\end{proof}

Theorem \ref{equiv_thm} follows from these three lemmata; the maps $\tau_n$ assemble into an equivalence of pro-ring spectra.  It is worth pointing out that these methods extend to give an equivalence $Ad(M \times_G EG)^{-TBG} \simeq M^{hG}$ of pro-ring spectra for any $G$-monoid $M$.


\section{Hochschild cohomology of $C^*(BG)$} \label{hoch_BG_section}

The purpose of this section is to prove the equivalence of parts
(\ref{gs_item}) and (\ref{hoch_BG_item}) in Theorem \ref{hoch_thm}.
We begin with a cosimplicial description of the terms in the
prospectrum $LBG^{-TBG}$.  We use this to give an intermediate result
describing the homology of these terms.  This is then assembled into
the result using various limit arguments.

Because we have assumed that $G$ is connected, $BG$ is simply connected, and for $n$ sufficiently large, so too is $BG_n$.  This ensures that the spectral sequences that we employ will converge.

\subsection{A cosimplicial model for $Ad(EG_n)^{-TBG_n}$}

In this section, we construct a cosimplicial ring spectrum $Ad_n^\bullet$ with the property that
$$\Tot(Ad_n^\bullet) \simeq Ad(EG_n)^{-TBG_n}$$
The bulk of this section is adapted directly from \cite{cj, cohenatiyah}, so we will be brief except in instances where our construction differs substantially.

One can realize the free loop space of $BG$ as the totalization of the cosimplicial space $Map(S^1_{\bullet}, BG)$:
$$LBG = Map(S^1, BG) = Map(|S^1_{\bullet}|, BG) = \Tot (Map(S^1_{\bullet}, BG))$$
Here $S^1_{\bullet}$ is the simplicial set whose geometric realization is the circle; $S^1_{\bullet}$ has $k+1$ $k$-dimensional simplices.  Hence
$$Map(S^1_k, BG) = BG^{\times k+1}$$
The cofaces and codegeneracies are given by various diagonals and projections.

\begin{proposition} \label{cosimp_prop}

The space $Ad(EG_n)$ is homotopy equivalent to the totalization of the subcosimplicial space of $Map(S^1_{\bullet}, BG)$ whose $k^{\rm th}$ space is
$$BG_n \times BG^{\times k}$$

\end{proposition}

\begin{proof}

The subcosimplicial space described is carried via the equivalence 
$$\Tot(Map(S^1_{\bullet}, BG)) = LBG$$
homeomorphically to the subspace $L_nBG \subseteq LBG$ given by those loops whose basepoint lies in $BG_n \subseteq BG$.

Recall that $Ad(EG_n) = EG_n \times_G G$; as such $Ad(EG_n)$ is the fiber product:
$$\xymatrix{
Ad(EG_n) \ar[r]^{\subseteq} \ar[d] & Ad(EG) \ar[d] \\
BG_n \ar[r]^{\subseteq} & BG
}$$
Similarly, $L_nBG$ is the fiber product
$$\xymatrix{
L_nBG \ar[r]^{\subseteq} \ar[d] & LBG \ar[d] \\
BG_n \ar[r]^{\subseteq} & BG
}$$
Since the fibrations $LBG \to BG$ and $Ad(EG) \to BG$ are equivalent, these fiber squares imply that $L_nBG$ and $Ad(EG_n)$ are equivalent.

\end{proof}

We now desuspend this construction by the tangent bundle of $BG_n$.  For this we need the following construction.  Consider the composite map
$$\xymatrix{
BG_n \ar[r]^-{\Delta} & BG_n \times BG_n \ar[r]^{1 \times i} & BG_n \times BG
}$$
This is the $0^{\rm th}$ coface of the cosimplicial space which totalizes to $L_nBG$.  The pullback of $TBG_n \times 0$ to $BG_n$ via this map is once again $TBG_n$.  Thus we have an induced map
$$\mu_R: BG_n^{-TBG_n} \to BG_n^{-TBG_n} \smash BG_+$$
Making the same construction with $1 \times i$ replaced by $i \times 1$ defines a similar map 
$$\mu_L: BG_n^{-TBG_n} \to BG_+ \smash BG_n^{-TBG_n}$$

We describe these maps by the element-theoretic formulae
$$\begin{array}{cc}
\mu_L(u) = (y_L, v_L), & \mu_R(u) = (v_R, y_R)
\end{array}$$
Though this does not quite make sense as spectra do not have elements, we hope the meaning is clear.

\begin{definition}

For a group $G$ and an integer $n>0$, define a cosimplicial spectrum $Ad_n^{\bullet}$ whose $k^{\rm th}$ term is
$$Ad_n^k := BG_n^{-TBG_n} \smash BG^{\times k}_+$$
with coface and codegeneracy maps defined by the element-theoretic formulae
\begin{eqnarray*}
\delta_0(u; x_1, \dots, x_{k-1}) & = & (v_R; y_R,  x_1, \dots, x_{k-1}) \\
\delta_i(u; x_1, \dots, x_{k-1}) & = & (u; x_1, \dots, x_{i-1}, x_i, x_i, x_{i+1}, \dots, x_{k-1}) \\
                         &   & 1 \leq i \leq k-1 \\
\delta_k(u; x_1, \dots, x_{k-1}) & = & (v_L; x_1, \dots, x_{k-1}, y_L) \\
\sigma_i(u; x_1, \dots, x_{k+1}) & = & (u; x_1, \dots, x_i, x_{i+2}, \dots, x_{k+1}) \\
                         &   & 0 \leq i \leq k
\end{eqnarray*}

\end{definition}

Define a map
$$m_{k, l} : Ad_n^k \smash Ad_n^l \to Ad_n$$
by the composite
$$\xymatrix{
(BG_n^{-TBG_n} \smash BG^{\times k}_+) \smash (BG_n^{-TBG_n} \smash BG^{\times l}_+) \ar[r]^-{T}  \ar[dr]_-{m_{k, l}} & BG_n^{-TBG_n}  \smash BG_n^{-TBG_n}  \smash BG^{\times k+l}_+ \ar[d]^-{m \smash 1} \\
 & BG_n^{-TBG_n}  \smash BG^{\times k+l}_+
}$$
where $T$ switches factors, and $m$ is multiplication in the ring spectrum $BG_n^{-TBG_n}$.

After totalization, the maps $m_{k, l}$ define a multiplication
$$\Tot(Ad_n^\bullet) \smash \Tot(Ad_n^{\bullet}) \to \Tot(Ad_n^\bullet)$$
which makes $\Tot(Ad_n^\bullet)$ into a ring spectrum.

\begin{theorem} \label{cosimp_spectrum_thm}

There is an equivalence of ring spectra
$$\xymatrix{
Ad(EG_n)^{-TBG_n} \ar[r]^-{\simeq} & \Tot (Ad_n^{\bullet})
}$$

\end{theorem}

\begin{proof}

The equivalence of these spectra follows from Proposition \ref{cosimp_prop}.  The proof that the equivalence preserves ring multiplication is identical to the proof in \cite{cj} that $LM^{-TM}$ and $\Tot(\L_M^\bullet)$ are equivalent ring spectra.

\end{proof}

\subsection{The homology of $Ad(EG_n)^{-TBG_n}$}

%
%
%

Recall that for any space $X$, the singular cochain complex, $C^*(X)$, is a differential graded algebra via the cup product of cochains.  Using left and right multiplication, $C^*(X)$ becomes a $C^*(X)$-differential-graded bimodule.  Maps of spaces induce maps of differential graded algebras, so the maps
$$\xymatrix{
BG_n \ar[r]^-{i_n} & BG_{n+1} \ar[r]^-{i} & BG
}$$
make $C^*(BG_n)$, $C^*(BG_{n+1})$ and $C^*(BG)$ into $C^*(BG)$-bimodule algebras.  Further, the maps $i_n^*$ and $i^*$ are maps of bimodule algebras.  One may therefore form the Hochschild cohomology
$$HH^*(C^*(BG), C^*(BG_n))$$
which becomes a ring under the cup product of Hochschild cochains.  This allows us to describe the homology of individual terms of the prospectrum $LBG^{-TBG}$:

\begin{theorem} \label{eg_n_thm}

There is a ring isomorphism
$$HH^*(C^*(BG), C^*(BG_n)) \cong H_*(Ad(EG_n)^{-TBG_n})$$

\end{theorem}

\begin{proof}

Theorem \ref{cosimp_spectrum_thm} gives the following equivalence of chain complexes:
$$C_*(Ad(EG_n)^{-TBG_n}) \simeq \Tot(C_*(BG^{\times \bullet}_+ \smash BG_n^{-TBG_n}))$$
Using the Eilenberg-Zilber theorem and Atiyah duality, the righthand side is equivalent to the totalization of the cosimplicial chain complex
$$k \mapsto C_*(BG)^{\otimes k} \otimes C^*(BG_n)$$

Define a chain map $g_k: C_*(BG)^{\otimes k} \otimes C^*(BG_n) \to \Hom(C^*(BG)^{\otimes k}, C^*(BG_n))$ by adjunction and evaluation:
$$e_1 \otimes \cdots \otimes e_k \otimes f \longmapsto ((f_1 \otimes \cdots \otimes f_k) \mapsto f_1(e_1) \cdots f_k(e_k) \cdot f)$$
It is easy to verify that the collection $\{ g_k, k \geq 0 \}$ defines a cosimplicial map
$$g: C_*(BG)^{\otimes \bullet} \otimes C^*(BG_n) \to CH^\bullet(C^*(BG), C^*(BG_n))$$
The theorem follows if we show that $g$ induces a homology isomorphism upon totalization.

To see this, we notice that there are spectral sequences that compute the homology of the two terms in question:
$$\begin{array}{ccc}
E_1 := H_*(BG)^{\otimes \bullet} \otimes H^*(BG_n) & \implies & H_*(\Tot(C_*(BG)^{\otimes \bullet} \otimes C^*(BG_n))) \\
 & {\rm and} & \\
E_1' := CH^*(H^*(BG), H^*(BG_n)) & \implies & HH^*(C^*(BG), C^*(BG_n)) 
\end{array}$$
The cosimplicial chain map $g$ induces a map $g_*$ between the spectral sequences; we claim that $g_*:E_1 \to E_1'$ is an isomorphism.  In each cosimplicial degree $k$, the map
$$g_*: H_*(BG)^{\otimes k} \otimes H^*(BG_n) \to \Hom(H^*(BG)^{\otimes k}, H^*(BG_n))$$
is a graded isomorphism because $H^*(BG_n)$ is finite dimensional, and $H_*(BG)^{\otimes k}$ is finite dimensional in each degree.  Consequently $g_*$ is an isomorphism of spectral sequences; hence $g$ induces an isomorphism in homology after totalization.

The cosimplicial product structure on $Ad_n^\bullet$ is seen immediately to coincide with the cup product of Hochschild cochains.  Consequently, this is an isomorphism of rings.

\end{proof}

\subsection{Limit arguments} \label{limit_section}

Examine the direct system
$$BG_1 \to \dots \to BG_n \to BG_{n+1} \to \dots \to BG$$
Applying the (integral) singular chain and cochain complex functors produces direct and inverse systems of chain (resp. cochain) complexes.  Since $BG$ is given the weak (or limit) topology of the system, this allows us to identify the singular chain complex of $BG$:
$$C_*(BG) = \varinjlim C_*(BG_n)$$
Standard properties of limits and colimits then imply that
$$C^*(BG) = \varprojlim C^*(BG_n)$$

\begin{proposition} \label{ch_lim_prop}

There is an isomorphism of cochain complexes
$$CH^*(C^*(BG), C^*(BG)) \cong \varprojlim CH^*(C^*(BG), C^*(BG_n))$$

\end{proposition}

\begin{proof}

For a given differential graded algebra $A$, the Hochschild cochain functor $CH^*(A, \cdot)$ is covariant in the module variable for chain maps of differential graded modules over $A$.

Recall that
$$i_n^*: C^*(BG_{n+1}) \to C^*(BG_n)$$
is a chain map and map of $C^*(BG)$-modules.  Consequently the map induced by $i_n^*$
$$CH^*(C^*(BG), C^*(BG_{n+1})) \to CH^*(C^*(BG), C^*(BG_n))$$
is a chain map.  Therefore $\varprojlim CH^*(C^*(BG), C^*(BG_n))$ is also a chain complex.

We also know that
$$i^*: C^*(BG) \to C^*(BG_n)$$
is a chain map and map of $C^*(BG)$-modules.  So the maps
$$CH^*(C^*(BG), C^*(BG)) \to CH^*(C^*(BG), C^*(BG_n))$$
are chain maps.  Since they are coherent across the inverse system, they assemble into a chain map
$$CH^*(C^*(BG), C^*(BG)) \to \varprojlim CH^*(C^*(BG), C^*(BG_n))$$

Generally, if $Z$ is an abelian group and
$$X_0 \gets X_1 \gets \cdots $$
an inverse system of abelian groups, there is a canonical isomorphism (of groups)
$$Hom(Z, \varprojlim X_i) \cong \varprojlim Hom(Z, X_i)$$
Consequently the map induced by $i^*$ is an isomorphism
\begin{eqnarray*}
CH^k(C^*(BG), C^*(BG)) & \cong & CH^k(C^*(BG), \varprojlim C^*(BG_n)) \\
 & \cong & \varprojlim CH^k(C^*(BG), C^*(BG_n))
\end{eqnarray*}
for each $k$.  The previous comments imply that this isomorphism is one of chain complexes.
%
%
%
%

\end{proof}

Using a $\varprojlim^1$ argument and some topology, we may conclude the following homological analogue.

\begin{corollary} \label{hh_lim_cor}

There is an isomorphism of rings
$$HH^*(C^*(BG), C^*(BG)) \cong \varprojlim HH^*(C^*(BG), C^*(BG_n))$$

\end{corollary}

\begin{proof}

The tower
$$\dots \gets C^*(BG_n) \gets C^*(BG_{n+1}) \gets \cdots$$
satisfies the Mittag-Leffler condition; consequently, so does the tower
$$\dots \gets CH^*(C^*(BG), C^*(BG_n)) \gets CH^*(C^*(BG), C^*(BG_{n+1})) \gets \cdots$$
Using this fact and the previous proposition, we see that there is a short exact sequence
\begin{eqnarray*}
0 \to {\varprojlim}^1 HH^*(C^*(BG), C^*(BG_n)) & \to & HH^*(C^*(BG), C^*(BG)) \\
 & \to & \varprojlim HH^*(C^*(BG), C^*(BG_n)) \to 0
\end{eqnarray*}

Recall that we've shown
$$HH^*(C^*(BG), C^*(BG_n)) \cong H_*(Ad(EG_n)^{-TBG_n})$$
So the $\varprojlim^1$ term vanishes if we can show that maps
$$Ad(EG_n)^{-TBG_n} \gets Ad(EG_{n+1})^{-TBG_{n+1}}$$
satisfy the Mittag-Leffler condition in homology.  Since the Spanier-Whitehead dual of $Ad(EG_n)^{-TBG_n}$ is $Ad(EG_n)^{-ad}$, this is equivalent to showing that the inclusions
$$Ad(EG_n) \to Ad(EG_{n+1}) \leqno{(*)}$$
satisfy the Mittag-Leffler condition in cohomology.  By construction, the connectivity of the inclusions $BG_n \hookrightarrow BG$ increases with $n$; hence the same is true for the inclusions $Ad(EG_n) \hookrightarrow Ad(EG)$.  This implies that $(*)$ does in fact satisfy the Mittag-Leffler condition in cohomology.

Since each map in the tower of coefficients is a ring homomorphism (in fact, a $C^*(BG)$-bimodule algebra map), the resulting isomorphism is one of rings.

\end{proof}

\subsection{A proof of (\ref{gs_item}) $\iff$ (\ref{hoch_BG_item}) in Theorem \ref{hoch_thm}}

Recall that we define
$$H_*^{pro}(LBG^{-TBG}) := \varprojlim H_*(Ad(EG_n)^{-TBG_n})$$
Using Theorem \ref{eg_n_thm} and Corollary \ref{hh_lim_cor}, we therefore have
$$H_*^{pro}(LBG^{-TBG}) \cong \varprojlim HH^*(C^*(BG), C^*(BG_n)) \cong HH^*(C^*(BG), C^*(BG))$$
The ring structure on the lefthand side is defined to be the inverse limit of the ring structures on $H_*(Ad(EG_n)^{-TBG_n})$.  We have just shown the same to be true for the righthand side; hence this isomorphism is one of rings.

\section{Spanier-Whitehead duality}

In this section, we show the isomorphism between the rings in parts
(\ref{gs_item}) and (\ref{fusion_item}) of Theorem \ref{hoch_thm}.  

Since $LBG \simeq Ad(EG) = \varinjlim \, Ad(EG_n)$, there is an exact
sequence
$$0 \to {\varprojlim}^1 H^*(Ad(EG_n)^{-ad}) \to H^*(LBG^{-ad}) \to
\varprojlim \,
H^*(Ad(EG_n)^{-ad}) \to 0$$
Using the same arguments as in section \ref{limit_section}, we see
that the ${\varprojlim}^1$ term vanishes.

In \cite{kate}, the first author has shown that the spectra
$Ad(EG_n)^{-TBG_n}$ and $Ad(EG_n)^{-ad}$ are Spanier-Whitehead dual.
Since these are finite spectra, we may conclude that
$$H_*(Ad(EG_n)^{-TBG_n}) \cong H^{-*}(Ad(EG_n)^{-ad})$$
Moreover, since Spanier-Whitehead duality carries the product on
$Ad(EG_n)^{-TBG_n}$ to the coproduct on $Ad(EG_n)^{-ad}$, this
isomorphism is one of rings.  Therefore there is a ring isomorphism
$$H_*^{pro}(LBG^{-TBG}) := \varprojlim \, H_*(Ad(EG_n)^{-TBG_n}) \cong
\varprojlim \, H^{-*}(Ad(EG_n)^{-ad}) \cong  H^{-*}(LBG^{-ad})$$

\section{Hochschild cohomology of $C_*(G)$ and Koszul duality} \label{hoch_G_section}

\subsection{The bar and cobar constructions}

We recall the bar and cobar constructions for differential graded
(co-)algebras.  To begin, let $R$ be an connected, augmented, associative dga over a field $k$.  Recall that for a right $R$-module $M$, and a left $R$-module $N$, the \emph{two-sided bar construction} $B(M, R, N)$ is the realization of the simplicial chain complex $B_\bullet(M, R, N)$, given by
$$\begin{array}{rl}
B_n(M, R, N) = M \otimes R^{\otimes n} \otimes N, & n \in \N
\end{array}$$
whose faces are given by multiplication in $R$ and the module
structure on $M$ and $N$ (and degeneracies are given by insertion of a
unit).  Recall, further, that $B(R) := B(k, R, k)$, the classic bar
construction on $R$, is a differential graded coalgebra, and $B(M, R, k)$ and $B(k, R, N)$ are, respectively, right and left comodules for $B(R)$.  

Dually, for a supplemented, coassociative coalgebra $S$ and right and left comodules $P$ and $Q$ for $S$, the \emph{two-sided cobar construction} $\Omega(P, S, Q)$ is the totalization of the cosimplicial chain complex
$$\begin{array}{rl}
\Omega^n(P, S, Q) = P \otimes S^{\otimes n} \otimes Q, & n \in \N
\end{array}$$
whose cofaces are given by comultiplication in $S$ and the comodule
structure on $P$ and $Q$, and whose codegeneracies come from the
counit in $S$.  Write $\Omega(S) := \Omega(k, S, k)$; this is a
differential graded algebra.

A relationship between these two constructions is as follows.  Let $S$ be a differential graded coalgebra over a field $k$ which is
finite dimensional in each degree.  Then the dual $S^\vee = \Hom(S,
k)$ is a differential graded algebra, and there is an isomorphism of
differential graded coalgebras:
$$B(S^\vee) \cong (\Omega(S))^\vee \leqno{(*)}$$

\subsection{Koszul duality}

To our knowledge, there are at least two approaches to proving that
Hochschild cohomology is insensitive to Koszul duality, using
\cite{fmt} and \cite{po}.  We recall these results.

A supplemented coalgebra $S = \overline{S} \oplus k$ is said to be
\emph{conilpotent} if, for every $x \in \overline{S}$, there is an $n$
so that the $n^{\rm th}$ iterated reduced comultiplication vanishes on
$x$. In \cite{fmt}, Felix, Menichi, and Thomas proved that if $S$ is
locally conilpotent, non-negatively graded, and finitely generated in
each degree, then there is an isomorphism of Gerstenhaber algebras
$$HH^*(\Omega S, \Omega S) \cong HH^*(S^\vee, S^\vee)$$
This was realized via a chain map
$$CH^*(\overline{\Omega} S, \overline{\Omega} S) \to CH^*(S^\vee, S^\vee)$$
Here $\overline{\Omega} S$ is the reduced cobar construction, which is
equivalent to $\Omega S$.

Dually, let $R$ be a differential graded algebra, and write $R^!$ for the Koszul dual dga of $R$.  That is, $R^!$ is the linear dual of $B(R)$:
$$R^! = (B(R))^\vee = \Hom(B(R), k)$$
In \cite{po}, Hu gave a proof that there is an equivalence of chain complexes
$$CH^*(R, R) \simeq CH^*(R^!, R^!)$$
assuming that $H_*(R^!)$ is a finite-dimensional $k$-vector space.
Though not explicitly stated, it is does follow from the
proof given there that this induces a ring isomorphism in Hochschild
cohomology (we include a sketch below).  These two results are clearly
related via the isomorphism $(*)$.

\begin{proposition}

The equivalence $CH^*(R, R) \simeq CH^*(R^!, R^!)$ of \cite{po}
induces a ring isomorphism
$$HH^*(R, R) \cong HH^*(R^!, R^!)$$

\end{proposition}

\begin{proof}

We summarize the essential points of the proof given in \cite{po} in order to show that this isomorphism is one of rings.  Hu considers the bicosimplicial object
$$X^{\bullet, \bullet} := \Hom_{R \otimes R^{op}}(B_\bullet(R, R, R), \Omega^\bullet(B(R, R, k), B(k, R, k), B(k, R, R)))$$
Recall that if $R$ is connected, there is an equivalence $R \to
\Omega(B(R))$ of differential graded algebras.  Further, there are
$R$-module equivalences $B(R, R, k) \to k \leftarrow B(k, R, R)$, so
we have an equivalence of $R \otimes R^{op}$-modules
$$R \to \Omega(B(R)) \leftarrow \Omega^\bullet(B(R, R, k), B(k, R, k), B(k, R, R)))$$
Furthermore, $B_\bullet(R, R, R)  \to R$ is a free $R\otimes
R^{op}$-resolution (over $B_\bullet(k, R, k)$).  Therefore
$X^{\bullet, \bullet}$ computes the Hochschild cohomology of $R$:
$$HH^*(R, R) = \RHom_{R \otimes R^{op}}(R, R) = H_*(X^{\bullet, \bullet})$$

Moreover, by $R\otimes R^{op}$-freeness, there is an isomorphism
\begin{eqnarray*}
X^{\bullet, \bullet} & = & \Hom_k(B_\bullet(k, R, k), \Omega^\bullet(B(R, R, k), B(k, R, k), B(k, R, R))) \\
 & = & \Hom_k(B_\bullet(R), \Omega^\bullet(B(R, R, k), B(R), B(k, R, R)))
\end{eqnarray*}
Using the equivalences $B(R, R, k) \simeq k \simeq B(k, R, R)$, we see
that this complex is equivalent to
$$\Hom_k(B_\bullet(R), \Omega^\bullet(k, B(R), k))$$
which is, in turn, isomorphic to
$$\Hom_{B(R) \otimes B(R)^{op}-comod}(B_\bullet(R), \Omega^\bullet(B(R), B(R), B(R)))  \leqno{(**)}$$
since $\Omega^\bullet(B(R), B(R), B(R)))$ is a cofree $B(R) \otimes B(R)^{op}$-comodule on $\Omega^\bullet(k, B(R), k)$.  Using the homological finiteness of $R^!$, we notice that $(**)$ computes
$$HH^*(R^!, R^!) = \RHom_{R^! \otimes {R^!}^{op}}(R^!, R^!) = \RHom_{B(R) \otimes B(R)^{op}-comod} (B(R), B(R))$$
because $\Omega^\bullet(B(R), B(R), B(R))$ is a cofree resolution of $B(R)$ in the category of $B(R) \otimes B(R)^{op}$-comodules.

Now, the Gerstenhaber cup product in $HH^*(A, A)$ can be identified with the Yoneda (composition) product in $\RHom_{A \otimes A^{op}}(A, A)$.  The isomorphism
$$HH^*(R, R) = H_*(X^{\bullet, \bullet}) \cong HH^*(R^!, R^!)$$
given above comes from the quasi-isomorphism
$$X^{\bullet, \bullet} \simeq \Hom_{B(R) \otimes B(R)^{op}-comod}(B_\bullet(R), \Omega^\bullet(B(R), B(R), B(R)))$$
which preserves the composition in each of these $\Hom$-complexes.

\end{proof}

\subsection{Application to $C^*(BG)$}

We will apply these results in the case at hand, using the coalgebra $S =
C_*(BG)$ or dually $R = S^\vee = C^*(BG)$.

It is well known (using the Eilenberg-Moore spectral sequence, for instance) that there is a homotopy equivalence of dga's
$$C_*(G) \simeq C_*(\Omega BG) \simeq \Omega(C_*(BG)) \simeq (C^*(BG))^!$$ 
and our assumption that $G$ is compact Lie ensures that the homology
of the Koszul dual is finite.  Using Hu's theorem, we conclude that
there is an ring isomorphism
$$HH^*(C^*(BG), C^*(BG)) \cong HH^*(C_*(G), C_*(G))$$

It is unclear whether we may employ \cite{fmt}to give an alternate
proof and strengthen this isomorphism to one of Gerstenhaber algebras.
For if we use the singular cochain complex, $S = C_*(BG)$ is far from
finite dimensional in each degree.  It may be possible to construct a
quasi-isomorphic coalgebra $S' \simeq S$ which satisfies the
assumptions of Felix-Menichi-Thomas' theorem (another of their results
implies that the Gerstenhaber structure of Hochschild cohomology is
preserved by quasi-isomorphism of dga's).  The simple
connectivity of $BG$ and local finiteness of $H_*(BG)$ suggest that
one may be able to find a locally finite simplicial set $Y_\bullet$
whose geometric realization is homotopy equivalent to $BG$.  Then $S'$
could be taken to be the simplicial chain complex of $Y_\bullet$.  But
we do not know a construction of such a simplicial set $Y_\bullet$.

\section{Relationship to string topology constructions}

We have already seen that several of our results have interpretations
in terms of string topology: in Theorem \ref{equiv_thm}, $LBG^{-TBG}$
arises from the string topology of $BG$, and the results of section \ref{hoch_BG_section} are analogues of the Cohen-Jones string topology theorem that for a simply connected manifold $M$, $\mathbb{H}_*(LM) \cong H^*(C^*(M), C^*(M))$ as graded algebras \cite{cj}.  In this section we will give an interpretation of the co-ring spectrum $LBG^{-ad}$ in terms of string topology, using string topology constructions for stacks.

In \cite{cg}, Cohen and Godin defined a non-counital Frobenius algebra structure on $h_*(LM)$, with multiplication given by the Chas-Sullivan product.  In \cite{lux}, Lupercio, Uribe and Xicot\'{e}ncatl extended the Chas-Sullivan construction to loop orbifolds.  Using this, a localization principle allowed them in \cite{lux_symm} to define an associative multiplication on $H_*(\Lambda [X^n/ \Sigma_n])$, the homology of the inertia orbifold of a symmetric product.  They then showed that this multiplication is Poincar\'{e} dual to a \emph{virtual intersection product} on $H^*(\Lambda [X^n/ \Sigma_n])$, which, with coauthors Gonz\'{a}lez and Segovia in \cite{glsux}, was identified with $H^*_{CR}(T^* [X^n/ \Sigma_n])$, the Chen-Ruan cohomology of the cotangent bundle of $[X^n/ \Sigma_n]$.  This product is part of a Frobenius algebra structure in Chen-Ruan cohomology.

Behrend, Ginot, Noohi, and Xu (BGNX) gave similar constructions in \cite{bgnxfrob, bgnxinertia}, where they define a Frobenius algebra structure on $H_*(\Lambda \mathfrak{X})$, the homology of the inertia stack of an oriented differentiable stack $\mathfrak{X}.$  Unlike the Frobenius algebra in Chen-Ruan cohomology, this structure is not necessarily unital nor counital.  In this structure, the multiplication is given by a stacky version of the Chas-Sullivan product, and the coproduct is given by a stacky version of the Cohen-Godin coproduct.

In the case that $\mathfrak{X} = [\ast /G ]$, the classifying stack of a compact Lie group $G$, the inertia stack $\Lambda \mathfrak{X}$ is the quotient stack $[G/G]$ where $G$ acts on itself by conjugation.  Then
$$ H_*(\Lambda \mathfrak{X} ) = H_*([G/G]) = H_*(Ad(EG)) \cong H_*(LBG) $$
so it is natural to ask whether the ``inertia Frobenius algebra" studied in \cite{glsux, bgnxfrob} is related to the co-ring spectrum  $LBG^{-ad}$.  The relationship is clearest when we consider instead the Frobenius algebra structure on $H^*(\Lambda [\ast/G]),$ induced via the universal coefficient theorem as in \cite{bgnxfrob}.  The following theorem says that the product (defined by BGNX) on $H^*([G/G])$, and hence the coproduct on $H_*([G/G]),$ are induced by the coproduct on the $LBG^{-ad}$ from Proposition~\ref{coring_prop}.

\begin{proposition}
The product on the inertia Frobenius algebra $H^*([G/G])$ is equal to the product on $H^*(LBG)$ induced from the co-ring spectrum structure on $LBG^{-ad}.$
\end{proposition} 

\begin{proof}
From Lemma 5.1 of \cite{bgnxfrob}, the product is given by
$$H^{i+j}([G \times G / G \times G]) \stackrel{\Delta^*}{\to} H^{i+j}([G \times G / G]) \stackrel{m_!}{\to} H^{i+j-d}([G / G])$$
where $d$ is the dimension of $G.$  Translating this product to homotopy orbit spaces gives:
$$H^{i+j}(G_{hG} \times G_{hG}) \stackrel{\Delta^*}{\to} H^{i+j}((G \times G)_{hG}) \stackrel{m_!}{\to} H^{i+j-d}(G_{hG}) $$
%
%
which is clearly the product given by applying $H^*$ and Thom isomorphisms to the coproduct on $LBG^{-ad}.$ 

\end{proof}

\begin{proposition}
There is a non-unital ring spectrum structure on $LBG^{+ad}$ which realizes the coproduct on $H^*([G/G])$ defined in \cite{bgnxfrob}.
\end{proposition}

\begin{remark}
However, BGNX have shown that this coproduct is trivial on $H^*([G/G];\mathbb{R}).$  It is likely to be nontrivial in any cohomology theory which detects the $G$-transfer map $\Sigma^{\infty} BG^{ad} \to S^0$ (such as orthogonal $K$-theory when $G=S^1$).
\end{remark}

\begin{proof}
The diagonal embedding $G \into G\times G$ induces a relative transfer map
$$ \tau_G^{G\times G}: (S^{\g \times \g} \wedge (G \times G)_+)_{hG\times G} \to (S^{\g}\wedge (G\times G)_+)_{hG}. $$
The lefthand side is equivalent to $(S^\g \wedge G_+)_{hG} \wedge (S^\g \wedge G_+)_{hG}.$  Group multiplication in $G$ induces
$$m: (S^{\g}\wedge (G\times G)_+)_{hG} \to (S^{\g}\wedge G_+)_{hG}$$
since it is $G$-equivariant.  Hence we can define the multiplication on 
$$LBG^{ad} = (S^\g \wedge G_+)_{hG} $$ 
to be
$$m \circ  \tau_G^{G\times G}: LBG^{ad} \wedge LBG^{ad} \to LBG^{ad}.$$
This product is the same as the ring structure on $LBG^{ad}$ described in \cite{goperads} coming from the first term of the transfer operad $G_{bG}$.  It is associative but not unital.  Applying cohomology and Thom isomorphisms yields 
$$H^i([G/G]) \stackrel{m^*}{\to} H^i([G \times G/G]) \stackrel{\tau^*}{\to} H^{i-d}([G \times G/G\times G])
\cong \bigoplus_{r+s = i-d} H^r([G/G]) \otimes H^s([G/G]) $$
which is the same as the description of the coproduct in Lemma 5.1 of \cite{bgnxfrob}.

\end{proof}

\bibliography{biblio}

\end{document}